\documentclass[11pt,reqno]{amsart}

\usepackage{latexsym}
\usepackage{amsfonts}
\usepackage{amsmath}
\usepackage{url}
\usepackage{nicefrac}
\usepackage{graphicx}
\usepackage{color}
\usepackage[colorlinks=true,citecolor=black,urlcolor=blue]{hyperref}

\usepackage{amssymb,amsmath}

\def\bl{\begin{lemma}}
\def\el{\end{lemma}}
\def\bth{\begin{theorem}}
\def\eth{\end{theorem}}
\def\bc{\begin{corollary}}
\def\ec{\end{corollary}}
\def\bcj{\begin{conjecture}}
\def\ecj{\end{conjecture}}
\def\bpr{\begin{proposition}}
\def\epr{\end{proposition}}
\def\bde{\begin{definition}}
\def\ede{\end{definition}}
\def\E{\mathbb{E}}

\newcommand{\be}{\begin{eqnarray}}
\newcommand{\ee}{\end{eqnarray}}

\newcommand{\R}{{\mathbb R}}

\newcommand{\Z}{{\mathbb Z}}

\renewcommand{\and}{\hbox{ {\rm and} }}

\newtheorem{theorem}{Theorem}
\newtheorem{definition}{Definition}[section]
\newtheorem{lemma}[theorem]{Lemma}

\newtheorem{corollary}[theorem]{Corollary}
\newtheorem{proposition}[theorem]{Proposition}
\newtheorem{conjecture}[theorem]{Conjecture}
\newtheorem{quest}[theorem]{Question}

\theoremstyle{definition}
\numberwithin{equation}{section}

\begin{document}
\title{Nonamenable Liouville Graphs}
\author{Itai Benjamini and Gady Kozma }

\begin{abstract}
Add to each level of binary tree edges to make the induced graph on the level a uniform expander.
It is shown that such a graph admits no non-constant bounded harmonic functions.
\end{abstract}

\maketitle

\section{Introduction }
We describe an example of a bounded degree graph with a positive Cheeger constant ({\it nonamenable})
which is {\it Liouville}, that is admits no non constant bounded harmonic functions.
A much more complicated example of bounded geometry simply connected Riemannian manifold with these properties
was constructed in \cite{BC}.
Nonamenable Cayley graphs are not Liouville  \cite{KV}, since the question ``does nonamenability implies non Liouville
for general bounded degree graphs?'' keeps coming up and the example below is transparent (unlike \cite{BC}) , we decided to write it down.
The last section contains related conjectures.

Recall that a graph is nonamenable if there is $C > 0$,
$$
\inf_{S}  { |\partial S| \over |S|} > C,
$$
for every finite nonempty set of vertices $S$, where $\partial S$ denotes the neighbours of $S$, (vertices not in $S$
with a neighbour in $S$).
\medskip

The basic idea is to start with a binary tree and add edges to it, in order
``to collapse the different directions simple random walk can escape to infinity''.
This will be done using {\it expanders}, which is a family of finite $d$-regular graphs, $G_n$, growing in size, for which the
isoperimetric condition above holds for any set of at most half the size of $G_n$, for all $n$.
For background on expanders see e.g. \cite{HLW}.

The example might be relevant to some modelling in genetics.

\section{The example}

\noindent
{\it Example:} On the vertices of each level of a binary tree place a $3$-regular expander.

\begin{proof}
Denote this graph by $G$.
Since the binary tree is nonamenable $G$ is too. For the Liouville
property, examine $L_n$, the $n^\textrm{th}$ level of the
tree. Consider it as a probability space (with the probability of each
point $2^{-n}$). For $x\in L_n$ and $y\in L_{n+1}$ define $p(x,y)$ be
the probability that random walk starting from $x$ will first hit
$L_{n+1}$ in $y$. The corresponding operator will be denoted by $T=T_n$
i.e.~for every $\phi:L_n\to\R$ let $T\phi:L_{n+1}\to\R$ be defined by
\[
(T\phi)(y)=2\sum_{x\in L_n}\phi(x)p(x,y)\qquad\forall y.
\]
Next, denote by $||T||_{p\to q}$ the norm of
$T$ as an operator from $L^p(L_n)\to L^q(L_{n+1})$. We immediately get
$||T||_{1\to 1}\le 1$ because $T$ is linear and sends probability
measures to probability measures (the $2$ in the definition of $T$
cancels out with the difference between the measures on $L_n$ and $L_{n+1}$).

Further, $||T||_{\infty\to\infty}$ is also $\le 1$. To see this write
\[
||T\phi||_\infty=2\max_{y\in L_{n+1}}\sum_{x\in L_n}\phi(x)p(x,y) \le
2||\phi||_\infty\max_{y\in L_{n+1}}\sum_{x\in L_n}p(x,y)
\]
so
\[
||T||_{\infty\to\infty}\le 2\max_{y\in L_{n+1}}\sum_xp(x,y).
\]
Now by time reversal, $\sum_x p(x,y)$ is the same as the expected
number of times a random walk starting from $y$ visits $L_n$ before
exiting the ball $B_n:=\cup_{k\le n}L_k$ --- here we need to define the
time the random walker exits $B_n$ only after the first step. This
expectation is easy to calculate. Denote the number of visits by
$V$. There is probability $\nicefrac{1}{6}$ that the first step of the
walker goes to $B_n$, which is necessary for $V\ne 0$. Afterwards each
time our walker reaches $L_n$ there is probability $\nicefrac{1}{3}$
for it to exit $B_n$. Hence the number of visits after first entering
$B_n$ is a geometric variable with mean 3, and we get $\E
V=\nicefrac{1}{6}\cdot 3=\nicefrac{1}{2}$ so
$||T||_{\infty\to\infty}\le 1$. By the Riesz-Thorin interpolation
theorem we get
\[
||T||_{2\to 2}\le 1.
\]

Let now $u\in L_k$ be some vertex and denote $\mu_n$
the harmonic measure $L_n$ starting at $u$ (assume $n\ge
k$). Let $f_n=2^n\mu_n-1$. In fact, $f_n$ is simply applying $T$
repeatedly to $2^k\mathbf{1}_u-1$ where $\mathbf{1}_u$ is a Kronecker
$\delta$ at $u$. Because $T$ sends probability measures to probability
measures we always get that $\sum_{x\in L_n}f(n)=0$. Because of the
probabilistic interpretation of $T$ we can examine the first step and
only then apply $T$. For the first step there is probability
$\nicefrac{1}{3}$ to go to $L_{n+1}$, probability $\nicefrac{1}{2}$ to
stay in $L_n$ and probability $\nicefrac{1}{6}$ to go back to
$L_{n-1}$. Denoting these operators (normalised to have norm 1) by
$S_1$, $S_2$ and $S_3$ we get
\[
(T_nf) =
\tfrac{1}{3}S_1f+\tfrac{1}{2}T_nS_2f+\tfrac{1}{6}T_nT_{n-1}S_3f.
\]
The same reasoning as above shows that $||S_i||_{2\to 2}\le 1$ and
further, because $S_2$ is the operator corresponding to a random walk
on an expander, and because the average of $f$ is zero we get
$||S_2f||_2\le (1-\lambda) f$ where $\lambda$ is the spectral gap of the
expander we put on the levels. In fact, let's define $\lambda$ as the
infimum of these spectral gaps over all levels so we do not need to
define $\lambda_n$ etc. All in all we get
\[
||Tf||_2 \le (1-\tfrac{1}{2}\lambda)||f||_2
\]
so $||f_n||_2\to 0$ as $n\to\infty$.

To finish the proof use that $||T||_{\infty\to\infty}\le 1$ so
$||f_n||_\infty\le C$. Hence we also get $||f_n||_1\to 0$ as
$n\to\infty$. To show that $G$ is Liouville, use Poisson's formula for
harmonic function $h$ with respect to level $n$ large.
\[
h(u) - h(v) = \sum_{w \in T_n} h(w) (\mu^u_n(w) - \mu^v_n(w) ) \leq \max_{ x \in G} h(x) ||\mu^u_n(w) - \mu^v(w)||_1 \underset{n}{\longrightarrow} 0,
\]
since $h$ is bounded. Here again $\mu^u_n $ and $\mu^v_n$ denote the harmonic measures on level $n$ starting at $u$ and $v$
respectively.
\end{proof}

\section{Conjectures and Questions}

What if we replace the binary tree by another nonamenable graph, on which we add expanders on all spheres around a fixed vertex?
We were able to show that there is such graph which is not Liouville.
The construction, very roughly, starts with an unbalanced tree (i.e.~a
binary tree where the father is connected to its left child by a
regular edge, but to its right child with a double edge). We then
utilised the fact that the harmonic measure is highly non-uniform on
the levels and added a graph which is an expander with respect to the
uniform measure but not with respect to the harmonic measure. We skip
all details.

\medskip
In the example balls around the root are uniform expanders, (when considering the induced finite graph on the vertices of the ball).
An {\it old conjecture} of the first author is that there is
no infinite Cayley graph with this property, moreover there is no infinite bounded degree connected graph in which all balls,
centered at any vertex, are uniformly expanders. Verify this for the example above. Here is a heuristic. 
A Cayley graph in which all balls are expanders is non amenable and will admit non constant bounded harmonic functions yet is somewhat similar to the example above and thus should be Liouville, contradiction.
\medskip

Maybe a stronger spectral requirement will imply existence of non constant bounded harmonic functions:

An  infinite $d$-regular graph  $G$  is called {\it Ramanujan}, if the spectral radius of the
Markov operator (acting on $l^2$ of the vertices) equals
$2 {\sqrt{d-1} \over d}$  which is what it is for the $d$-regular tree.
When $G$ is connected, this spectral radius can be expressed as the
limiting exponent of the probabilities of return of a random walk,
that is,  $ \lim p_n(v,v)^{1/n}$.

\begin{conjecture}
Infinite Ramanujan graphs are not Liouville.
\end{conjecture}

See \cite{GKN} for partial affirmative answer.

For general graphs,  Liouville property is unstable under quasi-isometries, \cite{T}. For a very simple example see \cite{BR}, \cite{BS}.
It will be useful to have,

\begin{conjecture}
A graph $G$, which is quasi-isometric to a nonamenable Cayley graph, is not Liouville.
\end{conjecture}

Conjecture 2 will be useful in pushing the heuristic regarding an attack on the old conjecture mentioned above for Cayley graphs by
adding more edges in a quasi isometric manner to get a graph closely  imitating the example above.
\medskip

Assume $\Z$ acts on $G$ by isometries and $H = G/ \Z$ is
Liouville. Moreover, assume that for some choice of fundamental domain
$F$, simple random walk on $G$ visits every translation of $F$ infinitely often a.s.

\begin{quest}
Is $G$ Liouville?
\end{quest}

In his proof of Gromov's theorem, Kleiner proved that every Cayley
graph admits a non-constant Lipschitz harmonic function,
\cite{Kl}. This cannot be true for a general graph because some
graphs admits no non-constant harmonic functions at all. For example,
the half-line, or consider a half ladder which admits an exponentially growing harmonic functions other then the constants. 
Uri Bader asked which general graphs  admit a non-constant Lipschitz harmonic function? In particular does every transient graph
or any graph with more than one end admit a non-constant Lipschitz harmonic function?

\subsection*{Acknowledgements}
IB would like to thank Tatiana Smirnova-Nagnibeda for inspiring
discussions. GK's research supported by the Israel Science Foundation.


\end{document}